\documentclass[11pt]{amsart}
\usepackage{geometry}                
\usepackage{graphicx}
\usepackage{amssymb}
\usepackage{epstopdf}
\newcommand{\diam}{\mbox{diam}}
\newcommand{\Pc}{\mathcal{P}}
\newcommand{\Fc}{\mathcal{F}}
\newcommand{\Cc}{\mathcal{C}}
\newcommand{\R}{\mathbb{R}}
\newcommand{\eps}{\varepsilon}
\newcommand{\Z}{\mathbb{Z}}

\title{Dobrushin and Steif metrics are equal}
\author{J. Armstrong-Goodall \and R.S.MacKay}
\address{Mathematics Institute, University of Warwick, Coventry CV4 7AL, UK}
\email{jacob.armstronggoodall@warwickgrad.net \and R.S.MacKay@warwick.ac.uk}
\date{\today}                                           

\begin{document}
\begin{abstract}
It is proved that two useful and apparently different metrics on the set of Borel probabilities on countable products of Polish spaces of bounded diameters are equal.  This paves the way for advances in their computation.
\end{abstract}

\maketitle
\section{Introduction}
It is well-known that standard metrics on spaces of multivariate probability distributions with many or countably infinite number of variables are of limited use.  For example, Liggett \cite{L} laments on p.70 that ``total variation convergence essentially never occurs for particle systems''.  Examples of other metrics that fail similarly are given in \cite{M1,M2}.

To rectify this, one of us \cite{M1} introduced a metric on multivariate probability distributions that {\em does} give convergence for many systems (the examples treated there were weakly dependent probabilistic cellular automata, but the same applies to particle systems).  It was based on ingredients from Dobrushin \cite{D} (following \cite{Va}), which give a type of weak convergence for such systems, but it appears that Dobrushin did not take the final step of metrising the weak topology (perhaps because it requires an assumption of bounded diameters).

It turned out, however, that Steif had proposed a metric that achieves the same goal many years before \cite{S}.  Its definition extends one of Ornstein \cite{OW} from the case of translation-invariant probabilities on $A^\Z$ for a finite set $A$ to general Borel probabilities on $A^\Z$.  \cite{M1} had dismissed extension of Ornstein's metric but had missed Steif's way of achieving it.

Superficially, the two metrics look different.  Yet in the Appendix to \cite{M2} it was shown that for finite spaces they are not only equivalent but equal.
This led to the conjecture that they are always equal.  

Dobrushin metric was defined on any countable product of Polish (complete separable metric) spaces with bounded diameters.
Steif's metric can easily be generalised to the same context.  In this paper it is proved that they are equal.

This result is significant because they are useful for proving and quantifying convergence of interacting particle systems and parameter-dependence of the stationary probabilities.  Their computation for explicit multivariate probabilities, however, is not easy (though see \cite{DM} for some successes).  It is helpful to have two alternative formulations of the same metric.

Before closing the introduction, a comment is appropriate on \cite{FH}.  It extends the weak convergence proof of \cite{Va} to allow a global component of interaction.  In our opinion, this is, however, already covered by the results of \cite{Va,D}, because the dependency matrix of \cite{Va,D} gains a contribution $\eps/N$ in each of the $N$ components, which still sums to only $\eps$.  Furthermore, the metric of \cite{FH} requires an artificial enumeration and weighting of the components, and they use the oscillation of a function rather than its Lipschitz constant. 


\section{Statement of result}
Let $S$ be a countable set.  For each $s \in S$, let $(X_s,d_s)$ be a Polish (complete separable metric) space.  Suppose $\sup_{s\in S}\diam(X_s) < \infty$.  Let $X = \prod_{s\in S} X_s$ with product topology.  Let $\Pc$ be the set of Borel probabilities on $X$.  For $\mu \in \Pc$ (or a signed Borel measure) and measurable $f:X \to \R$, denote the integral of $f$ with respect to $\mu$ by $\mu(f)$.
For $x \in X$ denote the component in $X_s$ by $x_s$.  
Extend $d_s: X_s \times X_s \to \R$ to a semi-metric $d_s: X \times X \to \R$ (denoted by the same symbol) defined by $d_s(x,y) = d_s(x_s,y_s)$ for all $x,y \in X$.

For $f: X \to \R$ and $s \in S$ define the {\em partial Lipschitz constant} $$\Delta_s(f) = \sup \frac{f(x)-f(y)}{d_s(x,y)}$$ over pairs $x\ne y \in X$ agreeing off $s$.  Let 
the {\em Dobrushin semi-norm} $$\|f\| = \sum_s \Delta_s(f) \in \R_+ \cup \{\infty\}.$$  Let the {\em Dobrushin smooth functions} $\Fc$ be the set of $f: X \to \R$ with $\|f\| < \infty$, and $\Cc$ be the constant functions $X \to \R$.  From these ingredients, \cite{M1} made the following
\vskip 1ex
\noindent{\bf Definition}: The {\em Dobrushin distance} between $\mu$ and $\nu \in \Pc$ is 
\begin{equation}
D(\mu,\nu) = \sup_{f \in \Fc\setminus \Cc} \frac{\mu(f)-\nu(f)}{\|f\|}.
\end{equation}

For $\mu,\nu \in \Pc$ let $M$ be the set of {\em joinings} of $\mu$ to $\nu$ (often called couplings), i.e.~the set of Borel probabilities $m$ on $X\times X$ whose marginals on the first and second factors are $\mu, \nu$, respectively.  Extending \cite{S} from the case where each $X_s$ was finite with discrete metric, define
\vskip 1ex
\noindent{\bf Definition}: The {\em Steif distance} between $\mu$ and $\nu \in \Pc$ is
\begin{equation}
\bar{d}(\mu,\nu) = \inf_{m \in M} \sup_{s\in S} m(d_s).
\end{equation}

It is not difficult to check (see \cite{M1} for $D$) that both $D$ and $\bar{d}$ are metrics on $\Pc$ and $\Pc$ is complete with respect to each.
\vskip 1ex
\noindent{\bf Theorem}:
$D=\bar{d}$

\section{Proof}
\label{sec:proof}

\noindent{\bf Proof}: Firstly, by homogeneity of degree one, 
$$D(\mu,\nu) = \sup_{f\in \Fc: \|f\|\le 1} (\mu-\nu)(f).$$ 
Let $E = \{e=(e_s)_{s\in S} \in \R_+^S: \sum_{s \in S} e_s \le 1\}$.  For $e \in E$, let $$c_e(x,y) = \sum_{s\in S} e_s d_s(x,y)$$ and
$$F_e = \{f:X\to \R: \forall x,y \in X, f(x)-f(y) \le c_e(x,y)\}.$$
Then $f \in \cup_{e\in E} F_e$ implies $\|f\|\le 1$ because if $f \in F_e$ then for all $s \in S$ and $x,y\in X$ agreeing on $s$, $f(x)-f(y) \le e_s d_s(x,y)$.  So $\Delta_s(f) \le e_s$, thus summing over $s\in S$, $\|f\|\le 1$.
Conversely, $\|f\| \le 1$ implies $f \in \cup_{e\in E} F_e$ because choose an enumeration of $S$ and change sequentially the components of $x$ to those of $y$ to obtain $f(x)-f(y) \le \sum_{s\in S} \Delta_s(f) d_s(x,y)$.  
But $\sum_{s\in S} \Delta_s(f) = \|f\| \le 1$ and $\Delta_s(f)\ge 0$, so $\Delta(f) \in E$, thus $f \in F_{\Delta(f)}$.
The supremum over $f \in \cup_{e\in E} F_e$ is the same as the supremum over $e\in E$ of the supremum over $f \in F_e$.
Thus, 
\begin{equation}
D(\mu,\nu) = \sup_{e\in E} \sup_{f\in F_e} (\mu-\nu)(f).
\label{eq:D2}
\end{equation}

Secondly, $\sup_{s\in S}m(d_s) = \sup_{e \in E}\sum_{s\in S} e_s m(d_s)$, because denote the lefthand side by $u$ and let $\eps>0$ then $\exists s' \in S$ with $m(d_{s'})\ge u-\eps$, so choose $e_{s'}=1$ and the rest of $e_s=0$ to get the righthand side at least $u-\eps$; conversely, $\sum_{s\in S}m(d_s) \le \sum_{s\in S} e_s u \le u$.  So 
$$\bar{d}(\mu,\nu) = \inf_{m\in M} \sup_{e\in E} \sum_{s\in S} e_s m(d_s).$$
This is at least $\sup_{e\in E} \inf_{m \in M} \sum_{s\in S} e_s m(d_s)$ because for any $e'\in E, m \in M$, $$\sup_{e\in E} \sum_{s \in S} e_s m(d_s) \ge \sum_{s\in S} e'_sm(d_s).$$  Taking the infimum over $m \in M$, $\inf_{m\in M} \sup_{e\in E} \sum_{s\in S} e_s m(d_s) \ge \inf_{m \in M} \sum_{s\in S} e'_s m(d_s)$.  Taking the supremum over $e' \in E$, $$\inf_{m\in M} \sup_{e\in E} \sum_{s\in S} e_s m(d_s) \ge \sup_{e' \in E} \inf_{m \in M} \sum_{s\in S} e'_s m(d_s).$$  
$\bar{d}(\mu,\nu)$ is also at most $\sup_{e\in E} \inf_{m \in M} \sum_{s\in S} e_s m(d_s)$ because for all $\eps>0$, $e \in E$, there exists $m'_{e,\eps}\in M$ such that 
$$\inf_{m \in M} \sum_{s \in S} e_s m(d_s) + \eps \ge \sum_{s\in S}e_s m'_{e,\eps}(d_s).$$  
So for all $\eps > 0$, 
$$\sup_{e \in E} \inf_{m \in M} \sum_{s\in S} e_s m(d_s) + \eps \ge \sup_{e\in E} \sum_{s\in S} e_{s} m'_{e,\eps}(d_s) \ge \inf_{m\in M} \sup_{e \in E} \sum_{s\in S} e_s m(d_s).$$  
This holds for all $\eps>0$ so the result $\bar{d}(\mu,\nu) \le \sup_{e\in E} \inf_{m \in M} \sum_{s\in S} e_s m(d_s)$ follows.
Combining the above two results, 
$$\bar{d}(\mu,\nu) = \sup_{e\in E} \inf_{m \in M} \sum_{s\in S} e_s m(d_s).$$
Thus, using linearity of integration and the definition of $c_e$, 
\begin{equation}
\bar{d}(\mu,\nu) = \sup_{e \in E} \inf_{m \in M} m(c_e).
\label{eq:db2}
\end{equation}

Thirdly, for all $e \in E$, $c_e$ is a semi-metric on $X$, so by Kantorovich-Rubinstein duality, e.g.~Theorem 5.10(i) of \cite{Vi},
\begin{equation}
\sup_{f \in F_e} (\mu-\nu)(f) = \inf_{m \in M} m(c_e).
\end{equation}
Taking the supremum over $e \in E$ and using equations (\ref{eq:D2},\ref{eq:db2}) yields the desired result:
$$D(\mu,\nu) = \bar{d}(\mu,\nu).$$ \qed

\section*{Acknowledgements}
We are grateful to Jeff Steif for bringing his metric to our attention and for correspondence on the topic.

\section*{Appendix: Clarification of Appendix of \cite{M2}}
We take the opportunity to clarify the Appendix to \cite{M2}.  For $c_e = \sum_{s\in S} e_s d_s$ and $\mu,\nu \in \Pc$, it addressed maximising $\mu(f)+\nu(g)$ over pairs of functions $f,g: X\to \R$ subject to $f(x)+g(y) \le c_e(x,y)$ for all $x,y \in X$.  To use linear programming results, the discussion there was restricted to the case of $X$ finite, but the analysis to follow here applies in full generality if maximum is replaced by supremum.

It was stated that ``For fixed $e$, the maximum is attained by $g=-f$, by the Kantorovich-Rubinstein theorem applied to cost function $\sum_{s\in S} e_s d_s(x_s,y_s)$."  What was intended to be cited is the third inequality in Theorem 5.10(i) of \cite{Vi}, combined with the statement there that one can impose $\psi$ to be $c$-convex, and an extension of the Particular Case 5.4 of \cite{Vi} to semi-metrics.  The extension was sketched in parentheses at the end of Particular Case 5.4 of \cite{Vi}, but we believe is missing a hypothesis, so we spell it out here.

Say $c: X\times X \to \R$ is a {\em semi-metric} if $\forall x,y,z \in X, c(x,z)\le c(x,y)+c(y,z)$ and $c(x,x)=0$.  Note that we do not require symmetry, nor non-negativity.
A function $\psi:X\to\R$ is called {\em $c$-convex} if there exists a function $\zeta:X\to \R$ such that $\psi(x) = \sup_y (\zeta(y)-c(x,y))$.
$\psi:X\to\R$ is called {\em 1-Lipschitz} (with respect to $c$) if for all $x,x' \in X$, $\psi(x)-\psi(x')\le c(x',x)$ (note that by reversing the roles of $x,x'$, it also implies $\psi(x)-\psi(x') \ge -c(x,x')$).
The {\em $c$-transform} of a function $\psi:X\to\R$ is the function $\psi^c(x) = \inf_y (\psi(y)+c(y,x))$.

\vskip 1ex
\noindent{\bf Proposition 1}: If $c$ is a semi-metric on $X$ and $\psi:X \to \R$, then the following are equivalent:
\begin{enumerate}
\item $\psi$ is $c$-convex;
\item $\psi$ is 1-Lipschitz;
\item $\psi^c = \psi$.
\end{enumerate}

\noindent{\bf Proof}: Suppose $\psi$ is $c$-convex.  Then for all $\eps>0$ and $x\in X$ there exists $y\in X$ such that $\psi(x) \le \zeta(y)-c(x,y)+\eps$.
Also, for all $x' \in X$, $\psi(x') \ge \zeta(y)-c(x',y)$, so $\psi(x)-\psi(x') \le c(x',y)-c(x,y)+\eps \le c(x',x)+\eps$ by the triangle inequality.  
So for all $x,x' \in X$, $\psi(x)-\psi(x')\le c(x',x)$, which is the definition of $\psi$ being 1-Lipschitz.

In the other direction, if $\psi$ is 1-Lipschitz, then for all $x,y \in X$, $\psi(y)-c(x,y)\le \psi(x)$ so $\sup_y (\psi(y)-c(x,y))\le \psi(x)$.  But if $y=x$ then $\psi(y)-c(x,y) = \psi(x)$, using $c(x,x)=0$.  Thus $\sup_y (\psi(y)-c(x,y)) = \psi(x)$, showing that $\psi$ is $c$-convex with $\zeta=\psi$.

Next, suppose $\psi$ is 1-Lipschitz.  Then for all $x,y \in X$, $\psi(x) \le \psi(y)+c(y,x)$.  So $\psi(x) \le \inf_y (\psi(y)+c(y,x))$.  The right hand side is the definition of $\psi^c(x)$.
Inserting $y=x$ and using $c(x,x)=0$ we see also that $\psi(x) \ge \psi^c(x)$.  So $\psi^c=\psi$.

Conversely, if $\psi=\psi^c$ then for all $x \in X$, $\psi(x)=\inf_y (\psi(y)+c(y,x))$.  Thus for all $\eps>0$ and $x \in X$, there exists $y$ with $\psi(y)+c(y,x) \le \psi(x)+\eps$.  Also, for all $x' \in X$, $\psi(x') \le \psi(y)+c(y,x')$.  So $\psi(x')-\psi(x) \le c(y,x')-c(y,x)+\eps \le c(x,x')+ \eps$ by the triangle inequality.  Hence $\psi(x')-\psi(x) \le c(x,x')$, so $\psi$ is 1-Lipschitz.
 \qed
\vskip 1ex

Then the third inequality of Theorem 5.10(i) of \cite{Vi} allows one to replace $g$ by $(-f)^c$
(the correspondence with his notation is $\phi = g, \psi = -f$).  
His statement there that one can impose that $\psi$ be $c$-convex and the above proposition allow one to replace $(-f)^c$ by $-f$.  Hence one can take $g=-f$, as claimed.

Note that Prop.5.8 of \cite{Vi} gives yet another equivalence to $c$-convexity: $\psi$ is $c$-convex iff $\psi^{cc} = \psi$.

We now give a simple proof that does not refer to the Kantorovich-Rubinstein theorem, using only the above Proposition 1.
\vskip 1ex
\noindent{\bf Proposition 2}: Given $\mu,\nu \in \Pc$ and semi-metric $c$ on $X$, the supremum of $\mu(f) + \nu(g)$ over pairs of functions $f,g:X \to \R$ subject to $f(x)+g(y)\le c(x,y)$ for all $x,y \in X$ is equal to the supremum over cases with $g=-f$.
\vskip 1ex
\noindent{\bf Proof}: Firstly, the constraint $f(x)+g(y) \le c(x,y)$ implies that for all $y$, $g(y) \le \inf_x (c(x,y)-f(x))$, which is the definition of $(-f)^c (y)$.  So $\nu(g) \le \nu((-f)^c)$ and for all $x,y \in X$, $f(x) + (-f)^c(y) \le c(x,y)$.  Thus replacing $g$ by $(-f)^c$ satisfies the constraint and does not decrease the objective function.

Secondly, by the same argument one can replace $f$ by $\tilde{f}=(-g)^c$.  But $-\tilde{f}$ is $c$-convex, because $(-g)^c(x)=\inf_y(-g(y)+c(y,x))$ iff $-(-g)^c (x) = \sup_y (g(y)-c(y,x))$.  Hence, by Proposition 1, $(-\tilde{f})^c = -\tilde{f}$.  So one can restrict to $f$ satisfying $(-f)^c=-f$.

Combining these, one can restrict $g$ to be $(-f)^c$ and we can restrict $f$ to be $(-(-f)^c)^c = f^c$, so we can restrict to $g=-f$.
\qed
\vskip 1ex

Finally, we remark that the final sentence of the Appendix to \cite{M2} perhaps left too much to the reader to do, but is fleshed out in section~\ref{sec:proof}.

\end{document}